# On the identification of the global reference set in data envelopment analysis

Mahmood Mehdiloozad[*]

*Department of Mathematics, College of Sciences, Shiraz University 71454, Shiraz, Iran*

S. Morteza Mirdehghan

*Department of Mathematics, College of Sciences, Shiraz University 71454, Shiraz, Iran*

Biresh K. Sahoo

*Xavier Institute of Management, Xavier University, Bhubaneswar 751 013, India*

Israfil Roshdi

*Department of Mathematics, Semnan Branch, Islamic Azad University, Semnan, Iran*

[*] **Corresponding author**: M. Mehdiloozad

Ph.D. student

Department of Mathematics

College of Sciences

Shiraz University

Golestan Street | Adabiat Crossroad | Shiraz 71454 | Iran

E-mail: m.mehdiloozad@gmail.com

Tel.: +98.9127431689

# On the identification of the global reference set in data envelopment analysis

## Abstract

It is well established that multiple reference sets may occur for a decision making unit (DMU) in the non-radial DEA (data envelopment analysis) setting. As our first contribution, we differentiate between three types of reference set. First, we introduce the notion of *unary reference set* (URS) corresponding to a given projection of an evaluated DMU. The URS includes efficient DMUs that are active in a specific convex combination producing the projection. Because of the occurrence of multiple URSs, we introduce the notion of *maximal reference set* (MRS) and define it as the union of all the URSs associated with the given projection. Since multiple projections may occur in non-radial DEA models, we further define the union of the MRSs associated with all the projections as unique *global reference set* (GRS) of the evaluated DMU. As the second contribution, we propose and substantiate a general linear programming (LP) based approach to identify the GRS. Since our approach makes the identification through the execution of a single primal-based LP model, it is computationally more efficient than the existing methods for its easy implementation in practical applications. Our last contribution is to measure returns to scale using a non-radial DEA model. This method effectively deals with the occurrence of multiple supporting hyperplanes arising either from multiplicity of projections or from non-full dimensionality of minimum face. Finally, an empirical analysis is conducted based on a real–life data set to demonstrate the ready applicability of our approach.

## 1. Introduction

Data envelopment analysis (DEA), introduced by Charnes, Cooper, and Rhodes (1978) and Charnes, Cooper, and Rhodes (1979) based on the seminal work of Farrell (1957), is a linear programming (LP) based method for measuring the relative efficiency of a homogeneous group of decision making units (DMUs) with multiple inputs and multiple outputs. Based on observed data and a set of postulates, DEA defines a reference technology set relative to which a DMU can be rated as *efficient* or *inefficient*. For an inefficient DMU, DEA recognizes a unique or multiple projection(s) on the efficient frontier of the technology set. Associated with each projection, it also identifies a set of observed efficient DMUs against which the under evaluation DMU is directly compared. Those efficient DMUs are called *reference DMUs*, and the corresponding set is referred to as a *reference set*.

The identification of *all* the possible reference DMUs for an inefficient unit is an important and interesting problem in DEA, on which we concentrate in this contribution by means of the non-radial range-adjusted model (RAM) of Cooper, Park, and Pastor (1999). This issue has received significant attention in the literature due to its wide range of potential applications in ranking (Jahanshahloo, Junior, Hosseinzadeh Lotfi, & Akbarian, 2007), benchmarking and target setting (Bergendahl, 1998; Camanho & Dyson, 1999), and measuring returns to scale (RTS) (Cooper, Seiford, & Tone, 2007; Krivonozhko, Førsund, & Lychev, 2014; Sueyoshi & Sekitani, 2007a; Sueyoshi & Sekitani, 2007b; Tone, 1996; Tone, 2005; Tone & Sahoo, 2006).

From a managerial point of view, the identification of all the reference DMUs is specifically important for two reasons. First, to improve the performance of an inefficient DMU, it may not be logical in practice to introduce an unobserved (virtual) projection as a benchmark. In such a situation, however, the identification provides the possibility to derive practical guidelines from benchmarking against the reference DMUs. Second, when some (but not all) reference DMUs are identified for an evaluated unit, the decision maker may be of the opinion that the identified DMUs are not appropriate benchmarks and may wish to have more options in choosing targets. In such a case, the identification allows him/her to incorporate the preference information into analysis so as to yield a projection with the most preferred (i) closeness (Tone, 2010), (ii) values of inputs and outputs, and (iii) shares of reference units in its formation.

The pioneer attempt to find all the reference DMUs in non-radial DEA models was made by Sueyoshi and Sekitani (2007b). Based on strong complementary slackness conditions (SCSCs) of linear programming, they proposed a primal–dual based method using the RAM model. The proposed method in their impressive study is very interesting as a theoretical idea. However, as Krivonozhko, Førsund, and Lychev (2012b) have argued, not only the computational burden of Sueyoshi and Sekitani's (2007b) approach is high, but it also seems that the basic matrices defined in their approach are likely to be ill-conditioned, leading to erroneous and unacceptable results even for medium-size problems. Furthermore, the economic interpretation of some constraints of their proposed model does not make sense. In a more recent and conscious attempt to overcome these difficulties, Krivonozhko et al. (2014) have proposed a primal–dual based procedure based on solving several LP problems. Using computational experiments, they showed that their proposed method works reliably and efficiently on real-life data sets and outperforms Sueyoshi and Sekitani's (2007b) approach.

It is worth noting that the studies conducted by Sueyoshi and Sekitani (2007) and Krivonozhko et al. (2014) correctly found all the observed DMUs on *minimum face* – a face of minimum dimension on which all the projections are located – as a *unique* reference set of a given DMU. On the other hand, both of these studies pointed out that the occurrence of *multiple* reference sets was possible. However, neither of them explicitly made a clear distinction between the uniquely-found reference set and other types of reference set for which multipleness may occur. This lack of discrimination creates an ambiguity about the uniqueness and, consequently, about the mathematical well-definedness of the definition of reference set.

Therefore, we were motivated to eliminate this ambiguity effectively. To do so, we have proposed three types of reference set sequentially, as our first contribution. Corresponding to a given projection, we first introduce the notion of *unary reference set* (URS) including efficient DMUs that are active in a specific convex combination producing this projection. Since multiple URSs (hereafter referred to as problem Type I) may occur, we introduce the notion of *maximal reference set* (MRS) and define it as the union of all the URSs associated with the given projection. Since multiple projections may occur in the RAM model, we further define the union of the MRSs associated with all the projections as unique *global reference set* (GRS) of the evaluated DMU. We have had an interesting finding: the convex hull of the GRS is equal to the minimum face. The benefits of the introduced three types of reference set (i.e., URS, MRS and GRS) are outlined below.

- The introduced concepts are all mathematically well-defined.
- The URS and MRS help demonstrate the occurrence of multiple reference sets associated with a single and multiple projection(s), respectively.
- While the multipleness may occur for the URS and MRS, the GRS presents a unique reference set that contains all the possible reference DMUs.

As our second contribution, we have proposed an LP model that identifies the GRS, and provides a projection in the relative interior of the minimum face. The proposed approach has several important features. First, it can effectively deal with the simultaneous occurrence of problems Type I and II. Second, this approach involves solving a single LP problem, which makes this approach computationally more efficient than the existing ones for its easy implementation in practical applications. Third, the computational efficiency of our approach is higher than that of the previous primal–dual ones, since it is developed based on the primal (envelopment) form that is computationally more efficient than the dual (multiplier) form (Cooper et al., 2007). Forth, since our proposed LP problem contains several upper-bounded variables, its computational efficiency can be enhanced by using the simplex algorithm adopted for solving the LP problems with upper-bounded variables, which is much more efficient than the ordinary simplex algorithm (Winston, 2003).

Fifth, our proposed approach is more general in the sense that it can be readily used without any change in both the 'additive model' (Charnes, Cooper, Golany, Seiford, & Stutz, 1985) and the 'BAM model' (Cooper, Pastor, J. T., Borras, Aparicio, & Pastor, D., 2011; Pastor, 1994; Pastor & Ruiz, 2007), because the difference between each of these two models and the RAM model lies only in the weights assigned to the input and output slacks in the objective function. With some minor changes, it can also be used in the 'RAM/BCC model' (Aida, Cooper, Pastor, & Sueyoshi, 1998), the 'DSBM model' of Jahanshahloo, Hosseinzadeh Lotfi, Mehdiloozad, and Roshdi (2012) and the 'GMDDF model' of Mehdiloozad, Sahoo, and Roshdi (2014). Furthermore, it can be easily implemented in any radial DEA model like the 'BCC model' of Banker, Charnes, and Cooper (1984), but with some minor changes. Finally, our proposed approach is free from the restricting assumption that the input–output data must be non-negative, so it can effectively deal with negative data. This can be very beneficial from a practical point of view since in many applications negative inputs or outputs could appear. See Pastor and Ruiz (2007) for various examples of applications with negative data.

The third contribution of this study is to measure the RTS in the non-radial DEA setting. As it is known, the concept of RTS is meaningful only when the relevant DMU lies on the frontier of the technology set. Hence, for an inefficient DMU, an efficient projection must be considered. In this case, the type and magnitude of the RTS is determined through the position(s) of the hyperplane(s) supporting the technology set at the projection used. The supporting hyperplane(s) passes/pass through the MRS associated with this projection and can be mathematically characterized via this MRS. Therefore, problem Type II causes the occurrence of multiple supporting hyperplanes (hereafter referred to as problem Type III), which makes the measurement of RTS difficult. Such a difficulty can be properly dealt with by using a relative interior point of the minimum face for the measurement of RTS. This is because the supporting hyperplane(s) binding at this point is/are characterized through the GRS, but not through a specific MRS. Nonetheless, the uniqueness of the characterized supporting hyperplane(s) cannot yet be guaranteed because the minimum face may not be a 'Full Dimensional Efficient Facet' (Olesen & Petersen, 1996; Olesen & Petersen, 2003).

To sum up, the difficulty raised by problem Type III in the measurement of RTS originates either from problem Type II or from the non-full dimensionality of the minimum face. To deal with this difficulty, we have developed a two–stage procedure for the measurement of RTS by exploiting the intensive study of Krivonozhko et al. (2014). In the first stage, we cope with the difficulty arising from problem Type II by finding a relative interior point of the minimum face via the LP problem proposed to identify the GRS. Then, for the obtained point[1], we use the indirect method of Banker, Cooper, Seiford, Thrall, and Zhu (2004) or the direct method of Førsund, Hjalmarsson, Krivonozhko, and Utkin (2007) to resolve the difficulty resulted from the non-full dimensionality of the minimum face. To demonstrate the ready applicability of our approach in empirical works, we have conducted an illustrative empirical analysis based on a real–life data set of 70 public schools in the United States.

The remainder of this paper unfolds as follows. Section 2 deals with the description of the technology followed by a brief review of the RAM model. Section 3 presents the main contribution of our study, where the three notions - URS, MRS, and GRS are introduced, and an LP model for the identification of the GRS is proposed. The model developed in this section is then used to develop a method for the measurement of RTS. Section 4 illustrates the application of our proposed approach with a numerical example, followed by an illustrative empirical application. Section 5 presents the summary of our work with some concluding remarks.

## 2 Background of the research

Throughout this paper we deal with $n$ observed DMUs; each uses $m$ inputs to produce $s$ outputs. For each DMU$_j$ ( $j \in J = \{1,...,n\}$ ), the input and output vectors are denoted by $\mathbf{x}_j = (x_{1j},...,x_{mj})^T$ $\in \mathbb{R}^m_{\geq 0}$ and $\mathbf{y}_j = (y_{1j},...,y_{sj})^T \in \mathbb{R}^s_{\geq 0}$, respectively. Superscript $T$ stands for a vector transpose. We have denoted the vectors and matrices in bold and have used $\mathbf{0}$ to show a vector with the value of 0 in every entry.

---

[1] Note that this point does not influence the RTS, since all the relative interior points of the minimum face have the same RTS (Krivonozhko, Førsund, & Lychev, 2012c).

## 2.1 Technology set

The technology set, $T$, is defined as the set of all feasible input–output combinations, i.e.,

$$T = \left\{ (\mathbf{x}, \mathbf{y}) \in \mathbb{R}^{m+s}_{\geq 0} \;\middle|\; \mathbf{x} \text{ can produce } \mathbf{y} \right\}. \tag{1}$$

Under the variable returns to scale (VRS) framework, the nonparametric DEA representation of $T$ can be set up as (Banker et al., 1984):

$$T^{DEA}_{VRS} = \left\{ (\mathbf{x}, \mathbf{y}) \in \mathbb{R}^{m+s}_{\geq 0} \;\middle|\; \sum_{j \in J} \lambda_j \mathbf{x}_j \leq \mathbf{x}, \; \sum_{j \in J} \lambda_j \mathbf{y}_j \geq \mathbf{y}, \; \sum_{j \in J} \lambda_j = 1, \; \lambda_j \geq 0, \; \forall j \in J \right\}. \tag{2}$$

**Definition 2.1.1** Let $H = \left\{ (\mathbf{x}, \mathbf{y}) \;\middle|\; \mathbf{u}^T (\mathbf{y} - \overline{\mathbf{y}}) - \mathbf{v}^T (\mathbf{x} - \overline{\mathbf{x}}) = 0 \right\}$ be a supporting hyperplane of $T^{DEA}_{VRS}$ at $(\overline{\mathbf{x}}, \overline{\mathbf{y}})$. Then, $H$ and its corresponding face, i.e., $F = H \cap T^{DEA}_{VRS}$, are called *strong* if and only if all the components of the coefficient vector $(\mathbf{u}, \mathbf{v})$ are positive.[2]

## 2.2 The RAM model

Considering DMU$_o$ ($o \in J$) to be the unit under evaluation, the RAM model (Cooper et al., 1999) is defined in reference to $T^{DEA}_{VRS}$ as

$$
\begin{aligned}
\rho_o = \min \quad & 1 - \frac{1}{m+s} \left\{ \sum_{i=1}^{m} \frac{s_i^-}{R_i^-} + \sum_{r=1}^{s} \frac{s_r^+}{R_r^+} \right\} \\
\text{s.t.} \quad & \sum_{j \in J} \lambda_j x_{ij} + s_i^- = x_{io}, \quad i = 1, \ldots, m, \\
& \sum_{j \in J} \lambda_j y_{rj} - s_r^+ = y_{ro}, \quad r = 1, \ldots, s, \\
& \sum_{j \in J} \lambda_j = 1, \\
& \lambda_j \geq 0, \; j \in J, \; s_i^-, s_r^+ \geq 0, \; \forall i, r,
\end{aligned}
\tag{3}
$$

---

[2] For more details, see Rockafellar (1970) and Davtalab-Olyaie, Roshdi, Jahanshahloo, and Asgharian (2014).

where $s_i^-$ $(\forall i)$ and $s_r^+$ $(\forall r)$ represent the excess of the $i$th input and the shortfall of the $r$th output, respectively. Here, $R_i^-$ $(\forall i)$ and $R_r^+$ $(\forall r)$ are the ranges defined by the lowest and highest observed values in the $i$th input and the $r$th output, respectively, i.e.,

$$R_i^- = \max_{j \in J}\left\{x_{ij}\right\} - \min_{j \in J}\left\{x_{ij}\right\}, \quad i = 1,...,m;$$
$$R_r^+ = \max_{j \in J}\left\{y_{rj}\right\} - \min_{j \in J}\left\{y_{rj}\right\}, \quad r = 1,...,s. \tag{4}$$

**Definition 2.2.1** DMU$_o$ is said to be *RAM-efficient* if and only if $\rho_o = 1$, i.e., all the slacks are zero at optimality in (3) (Brockett, Cooper, Golden, Rousseau, and Wang, 2004).

Let $\left(\boldsymbol{\lambda}^*, \mathbf{s}^{-*}, \mathbf{s}^{+*}\right)$ be an optimal solution to (3). Then, the projection of DMU$_o$ is defined by

$$P := \left(\hat{\mathbf{x}}_o, \hat{\mathbf{y}}_o\right) = \left(\mathbf{x}_o - \mathbf{s}^{-*}, \mathbf{y}_o + \mathbf{s}^{+*}\right) = \sum_{j \in J} \lambda_j \left(\mathbf{x}_j, \mathbf{y}_j\right). \tag{5}$$

It can then be easily proved that the projection $P$ is RAM-efficient.

## 3. Identification of the global reference set

### 3.1. The global reference set

In this subsection, we present some key definitions, concepts and results, which are all essential for the development of our proposed approach.

**Definition 3.1.1** Let $\left(\boldsymbol{\lambda}^*, \mathbf{s}^{-*}, \mathbf{s}^{+*}\right)$ be an optimal solution to (3) that is associated with a given projection $P$. We define the set of DMUs with positive $\lambda_j^*$ as the *unary reference set (URS)* for DMU$_o$ and denote it by $R_{oP}^*$ as

$$R_{oP}^* = \left\{\text{DMU}_j \mid \lambda_j^* > 0\right\}. \tag{6}$$

We refer to each member of $R_{oP}^*$ as a *reference DMU* of DMU$_o$. All the reference units of DMU$_o$ are RAM-efficient, and are located on a supporting hyperplane of $T_{VRS}^{DEA}$.

Since the projection $P$ may be expressed as multiple convex combinations of its associated reference DMUs, multiple optimal values may take place for the vector $\boldsymbol{\lambda}$, leading to the occurrence of multiple URSs (problem Type I). Under such an occurrence, the measurement of RTS via the approach of Tone (1996) may be problematic. For a detailed discussion on this issue, interested readers may refer to the illustrative Figures 1 and 2 in Krivonozhko et al. (2012c).

To deal with the occurrence of multiple URSs for a given projection $P$, we need to define a reference set containing all the possible URSs.

**Definition 3.1.2** We define the union of *all* the URSs associated with a given projection $P$ as the *maximal reference set (MRS)* for DMU$_o$ and denote it by $R_{oP}^M$ as

$$R_{oP}^M = \left\{ \mathrm{DMU}_j \,\middle|\, \lambda_j^* > 0 \text{ in some optimal solution of (3) associated with } P \right\}. \qquad (7)$$

Because the RAM model is non-radial in nature, it may produce multiple projections for DMU$_o$, resulting in the occurrence of multiple MRSs (problem Type II). The simultaneous occurrence of multiple URSs and multiple projections is illustrated with the help of an example in Section 3 in Sueyoshi and Sekitani (2007b).

To deal with the occurrence of multiple MRSs, we use the concept of *minimum face*[3] that was considered in detail by Sueyoshi and Sekitani (2007b) and Krivonozhko et al. (2014) from different sides. First, we formulate the set of all the optimal solutions of (3), $\Omega_o$, in the form of

$$
\Omega_o = \left\{ \left( \boldsymbol{\lambda}, \mathbf{s}^-, \mathbf{s}^+ \right) \,\middle|\, 
\begin{aligned}
& \sum_{j \in J_E} \lambda_j x_{ij} + s_i^- = x_{io}, \quad i = 1, \dots, m, \\
& \sum_{j \in J_E} \lambda_j y_{rj} - s_r^+ = y_{ro}, \quad r = 1, \dots, s, \\
& \sum_{j \in J_E} \lambda_j = 1, \\
& \sum_{i=1}^m \frac{s_i^-}{R_i^-} + \sum_{r=1}^s \frac{s_r^+}{R_r^+} = (m+s)(1-\rho_o), \\
& \lambda_j \geq 0, \; j \in J_E, \; s_i^-, s_r^+ \geq 0, \; \forall i, r
\end{aligned}
\right\}, \qquad (8)
$$

where $J_E$ is the index set of all RAM-efficient DMUs.

Then, the set of all the projections of DMU$_o$, referred to as a *projection set*, can be expressed as

$$\Lambda_o = \left\{ \left( \mathbf{x}_o - \mathbf{s}^-, \mathbf{y}_o + \mathbf{s}^+ \right) \,\middle|\, \left( \boldsymbol{\lambda}, \mathbf{s}^-, \mathbf{s}^+ \right) \in \Omega_o \right\}. \qquad (9)$$

---

[3] For a graphical illustration of the minimum face, see Figure 4 in Sueyoshi and Sekitani (2007b).

As demonstrated by Krivonozhko et al. (2014), there exists a face of minimum dimension, $\Gamma_o^{\min}$, which contains the projection set $\Lambda_o$. This face is referred to as minimum face and is, indeed, the intersection of all the faces of $T_{VRS}^{DEA}$ that contain $\Lambda_o$, i.e.,

$$\Gamma_o^{\min} = \bigcap_{\substack{F \text{ is a face of } T_{VRS}^{DEA} \\ \text{and } \Lambda_o \subseteq F}} F \,. \tag{10}$$

Now, we provide the following definition that considers the occurrence of multiple MRSs associated with multiple projections.

**Definition 3.1.3** We define the union of *all* the MRSs of DMU$_o$ as its *global reference set (GRS)* and denote it by $R_o^G$ as

$$R_o^G = \bigcup_{P \in \Lambda_o} R_{oP}^M \,. \tag{11}$$

**Lemma 3.1.1** The *convex hull* of the GRS, $conv\left(R_o^G\right)$, is a strong face of $T_{VRS}^{DEA}$.

See Appendix A for the proof.

**Theorem 3.1.1** The minimum face is equal to the convex hull of the GRS, i.e., $\Gamma_o^{\min} = conv\left(R_o^G\right)$.

See Appendix A for the proof.

This theorem reveals that the minimum face is spanned by the GRS. Specifically, we obtain the following two corollaries of Theorem 3.1.1.

**Corollary 3.1.1** An explicit representation of the minimum face is set up as

$$\Gamma_o^{\min} = \left\{ (\mathbf{x}, \mathbf{y}) \,\middle|\, \mathbf{x} = \sum_{j \in J_o^G} \lambda_j \mathbf{x}_j, \ \mathbf{y} = \sum_{j \in J_o^G} \lambda_j \mathbf{y}_j, \sum_{j \in J_o^G} \lambda_j = 1, \ \lambda_j \geq 0, \ \forall j \in J_o^G \right\}, \tag{12}$$

where $J_o^G \subseteq J_E$ is the index set of DMUs in $R_o^G$.

**Corollary 3.1.2** The minimum face is a polytope.[4]

---

[4] This result was proved in a different way by Krivonozhko et al. (2014).

The following lemma is a straightforward consequence of Definitions 3.1.2 and 3.1.3 and the expression (9).

**Lemma 3.1.2** $DMU_k \in R_o^G$ *if and only if* $\lambda_k > 0$ *in some* $\left(\boldsymbol{\lambda}, \mathbf{s}^-, \mathbf{s}^+\right) \in \Omega_o$. *Formally,*

$$R_o^G = \bigcup_{\left(\boldsymbol{\lambda}, \mathbf{s}^-, \mathbf{s}^+\right) \in \Omega_o} \left\{ DMU_j \,\middle|\, \lambda_j > 0 \right\}. \tag{13}$$

**Theorem 3.1.2** Let $\left(\boldsymbol{\lambda}', \mathbf{s}^{-\prime}, \mathbf{s}^{+\prime}\right)$ be an element of $\Omega_o$ such that $\boldsymbol{\lambda}'$ has the maximum number of positive components. Then,

$$R_o^G = \left\{ DMU_j \,\middle|\, \lambda_j' > 0 \right\}. \tag{14}$$

See Appendix A for the proof.

### 3.2. Identification of the global reference set

Consider the following homogeneous system of equations

$$\mathbf{Au} = \mathbf{0}, \ \mathbf{u} \geq \mathbf{0}, \tag{15}$$

where $\mathbf{A}$ is a matrix of order $p \times q$. Bertsimas and Tsitsiklis (1997) presented an LP problem to find a feasible solution vector $\mathbf{u} \in \mathbb{R}^q$ for (15) such that the number of its positive components is maximum (Exercise 3.27, pp. 136). In addition, they suggested the formulation of an LP problem to find a solution vector $\mathbf{u} \in \mathbb{R}^q$ with the maximum number of positive components for the non-homogeneous system of equations

$$\mathbf{Au} = \mathbf{d}, \ \mathbf{u} \geq \mathbf{0}, \tag{16}$$

where $\mathbf{d}$ is a given $p \times 1$ vector.

We extend the above-mentioned exercise in a way to help us develop an LP-based approach for the identification of the GRS.

**Lemma 3.2.1** Let $X$ be the set of feasible solutions of the following homogeneous system of equations

$$\mathbf{Au} + \mathbf{Bv} = \mathbf{0}, \ \mathbf{u} \geq \mathbf{0}, \ \mathbf{v} \geq \mathbf{0}, \tag{17}$$

where $\mathbf{A}$ and $\mathbf{B}$ are matrices of order $p \times q_1$ and $p \times q_2$, respectively. Further, let $\left(\mathbf{u}^*, \mathbf{w}^*, \mathbf{v}^*\right)$ be an optimal solution for the following LP problem:

$$
\begin{aligned}
\max \quad & \sum_{j=1}^{q_1} w_j \\
\text{s.t.} \quad & \sum_{j=1}^{q_1} \mathbf{a}_j \left(u_j + w_j\right) + \sum_{j=1}^{q_2} \mathbf{b}_j v_j = \mathbf{0}, \\
& 0 \leq w_j \leq 1, \ u_j \geq 0, \ j = 1, ..., q_1, \\
& v_j \geq 0, \ j = 1, ..., q_2,
\end{aligned}
\tag{18}
$$

where $\mathbf{a}_j$ and $\mathbf{b}_j$ denote the $j$th columns of the matrices $\mathbf{A}$ and $\mathbf{B}$, respectively. Then, $\left(\mathbf{u}^* + \mathbf{w}^*, \mathbf{v}^*\right)$ $\in X$, and $\mathbf{u}^* + \mathbf{w}^*$ has the maximum number of positive components.

See Appendix A for the proof.

**Theorem 3.2.1** Let $X$ be the non-empty set of feasible solutions of the following non-homogeneous system of equations

$$
\mathbf{A}\mathbf{u} + \mathbf{B}\mathbf{v} = \mathbf{d}, \ \mathbf{u} \geq \mathbf{0}, \ \mathbf{v} \geq \mathbf{0} .
\tag{19}
$$

Further, let $\left(u_j^*, w_j^*, j = 1, ..., q_1 + 1, v_j^*, j = 1, ..., q_2\right)$ be an optimal solution to the following LP problem:

$$
\begin{aligned}
\max \quad & \sum_{j=1}^{q_1+1} w_j \\
\text{s.t.} \quad & \sum_{j=1}^{q_1} \mathbf{a}_j \left(u_j + w_j\right) + \sum_{j=1}^{q_2} \mathbf{b}_j v_j - \mathbf{d}\left(u_{q_1+1} + w_{q_1+1}\right) = \mathbf{0}, \\
& 0 \leq w_j \leq 1, \ u_j \geq 0, \ j = 1, ..., q_1 + 1, \\
& v_j \geq 0, \ j = 1, ..., q_2.
\end{aligned}
\tag{20}
$$

Then, $\left(u_j' := \dfrac{u_j^* + w_j^*}{u_{q_1+1}^* + w_{q_1+1}^*}, \ j = 1, ..., q_1, \ v_j' := \dfrac{v_j^*}{u_{q_1+1}^* + w_{q_1+1}^*}, \ j = 1, ..., q_2 \right) \in X$, and $\mathbf{u}'$ has the maximum number of positive components.

See Appendix A for the proof.

According to Theorem 3.1.2, a way of identifying the GRS is to find an element, namely $\left(\boldsymbol{\lambda}', \mathbf{s}^{-\prime}, \mathbf{s}^{+\prime}\right)$, in $\Omega_o$ such that the number of positive components of $\boldsymbol{\lambda}'$ is maximum. On the basis of this finding, we use Theorem 3.2.1 to formulate the LP problem

$$\max \quad \sum_{j \in J_E} \mu_j + \mu_{j_{t+1}}$$

$$\text{s.t.} \quad \sum_{j \in J_E} \left( \lambda_j + \mu_j \right) x_{ij} + s_i^- - \left( \lambda_{j_{t+1}} + \mu_{j_{t+1}} \right) x_{io} \quad = 0, \qquad i = 1, ..., m,$$

$$\sum_{j \in J_E} \left( \lambda_j + \mu_j \right) y_{rj} - s_r^+ - \left( \lambda_{j_{t+1}} + \mu_{j_{t+1}} \right) y_{ro} \quad = 0, \qquad r = 1, ..., s,$$

$$\sum_{j \in J_E} \left( \lambda_j + \mu_j \right) \qquad - \left( \lambda_{j_{t+1}} + \mu_{j_{t+1}} \right) \qquad = 0, \qquad (21)$$

$$\sum_{i=1}^{m} \frac{s_i^-}{R_i^-} + \sum_{r=1}^{s} \frac{s_r^+}{R_r^+} \quad - \left( \lambda_{j_{t+1}} + \mu_{j_{t+1}} \right) (m+s)(1-\rho_o) = 0,$$

$$0 \leq \mu_{j_k} \leq 1, \ \lambda_{j_k} \geq 0, \ k = 1, ..., t+1,$$

$$s_i^- \geq 0, \ s_r^+ \geq 0, \ i = 1, ..., m, \ r = 1, ..., s,$$

where $t$ denotes the cardinality of $J_E$, i.e., $J_E = \left\{ j_1, ..., j_t \right\}$.

Let $\left( \lambda_{j_k}^*, \mu_{j_k}^*, s_i^{-*}, s_r^{+*}, \forall k, i, r \right)$ be an optimal solution to (21). Then, Theorem 3.2.1 follows that

$$\left( \lambda_{j_k}' = \frac{\lambda_{j_k}^* + \mu_{j_k}^*}{\lambda_{j_{t+1}}^* + \mu_{j_{t+1}}^*}, s_i^{-\prime} = \frac{s_i^{-*}}{\lambda_{j_{t+1}}^* + \mu_{j_{t+1}}^*}, s_r^{+\prime} = \frac{s_r^{+*}}{\lambda_{j_{t+1}}^* + \mu_{j_{t+1}}^*}, \forall k, i, r \right) \in \Omega_o \qquad (22)$$

and $\lambda'$ has the maximum number of positive components. Applying Theorem 3.1.2 to (22), we identify the GRS as follows:

$$R_o^G = \left\{ \text{DMU}_j \, \middle| \, \lambda_j' > 0 \text{ in (22)} \right\}. \qquad (23)$$

Having identified the GRS, the projection associated with (22) can also be obtained as

$$P_o^* = \left( \mathbf{x}_o^*, \mathbf{y}_o^* \right) = \left( \mathbf{x}_o - \mathbf{s}^{+\prime}, \mathbf{y}_o + \mathbf{s}^{-\prime} \right) = \sum_{j \in J_o^G} \lambda_j' \left( \mathbf{x}_j, \mathbf{y}_j \right). \qquad (24)$$

Since $P_o^*$ is expressed as a strict convex combination of the units in $R_o^G$, it is a relative interior point of $\Gamma_o^{\min}$.

In summary, the solutions to (21) determine the set of all the DMUs spanning the minimum face (i.e., the GRS) together with a relative interior point of this face (i.e., the projection $P_o^*$).

### 3.3. Properties of the proposed approach

Some useful properties of the proposed approach are presented below.

➢ *Computational efficiency*

Our approach for determining the GRS involves the execution of a single LP problem, which makes it computationally more efficient than the existing ones for its easy implementation in practical applications. Moreover, the computational efficiency of the proposed approach is higher than that of the previous primal–dual ones since it is based on the primal (envelopment) form which is computationally more efficient than the dual (multiplier) form (Cooper et al., 2007).

Furthermore, since our proposed LP problem contains several upper-bounded variables, its computational efficiency can be enhanced by using the simplex algorithm adopted for solving the LP problems with upper-bounded variables. This is because considering (21) as an LP problem with upper-bounded variables leads to a further reduction in its size. More precisely, the size of the basic matrices during the solution process becomes $(m+s+2)\times(m+s+2)$, which is greater than the size of the basic matrices in model (3) by 3.

➤ *RTS Measurement*

As is well known, the concept of RTS is meaningful only when the relevant DMU lies on the frontier of the technology set. Hence, for an inefficient DMU$_o$, an efficient projection must be considered. In this case, the type and magnitude of the RTS for DMU$_o$ is determined through the position(s) of the hyperplane(s) supporting $T_{VRS}^{DEA}$ at the projection used. The supporting hyperplane(s) passes/pass through the MRS associated with this projection and can be mathematically characterized via this MRS. This implies that the RTS measurement would be problematic under the occurrence of problem Type II, as it causes the occurrence of multiple MRSs and, consequently, the occurrence of multiple supporting hyperplanes (problem Type III). Note that overcoming the difficulty caused by problem Type II does not necessarily guarantee the uniqueness of the characterized supporting hyperplane(s). This is because the minimum face may not be a 'Full Dimensional Efficient Facet' – an efficient facet of dimension $m+s-1$ in the input–output space (Olesen & Petersen, 1996, 2003).

To sum up, the difficulty in the measurement of RTS arises either from Type II problem or from the non-full dimensionality of the minimum face. To overcome this difficulty, we propose a two–stage procedure based on the intensive study of Krivonozhko et al. (2014). The first stage is to determine a relative interior point of the minimum face by (24). This determination is based on the following two facts which show that the difficulty arising from problem Type II can be dealt with by using a relative interior point of the minimum face for the measurement of RTS.

(a) All relative interior points of the minimum face $\Gamma_o^{\min}$ operate under the same type of RTS (Krivonozhko et al, 2012c).

(b) By Theorem 3.1.1, any relative interior point of $\Gamma_o^{\min}$ can be expressed as a strict convex combination of all the reference DMUs in $R_o^G$. Therefore, each supporting hyperplane binding at this point is also binding at all the reference DMUs. Precisely, this supporting hyperplane is characterized through the GRS, but not through a specific MRS.

In the second stage, the difficulty resulted from the non-full dimensionality of the minimum face is dealt with by measuring the RTS of $P_o^*$ through the indirect method of Banker et al. (2004) or the direct method of Førsund et al. (2007).

➢ *Extension to other DEA models*

Our approach can readily be used without any change for the 'additive model' (Charnes et al., 1985) and the 'BAM model' (Cooper et al., 2011; Pastor, 1994; Pastor & Ruiz, 2007), because the difference between each of these two models and the RAM model lies only in the weights assigned to the input and output slacks in the objective function. With some minor changes, it can also be adopted for the 'RAM/BCC model' of Aida et al. (1998), the 'DSBM model' of Jahanshahloo et al. (2012) and the 'GMDDF model' of Mehdiloozad et al. (2014).

Furthermore, it can be implemented in any radial DEA model (e.g., the BCC model[5]). In this regard, let $\left(\theta^*, \boldsymbol{\lambda}^*, \mathbf{s}^{-*}, \mathbf{s}^{+*}\right)$ be an optimal solution to the input-oriented BCC model. Then, it is sufficient to replace $\mathbf{x}_o$ by $\theta^* \mathbf{x}_o$ in the input constrains of (8) and (21), and fix the sum of slacks at $\sum_{i=1}^m s_i^{-*} + \sum_{r=1}^s s_r^{+*}$ .

➢ *Extension to constant returns to scale case*

The assumption of VRS is maintained in our study. This is because when a data set contains some negative values, one may not be able to define an efficient frontier passing through the origin, as is assumed under constant returns to scale (CRS). Therefore, as argued by Silva Portela and Thanassoulis (2010), the assumption of CRS is untenable with negative data.

---

[5] The issue of identifying all possible reference DMUs of each inefficient unit using the BCC model was explored by Jahanshahloo, Shirzadi, and Mirdehghan (2008), Krivonozhko, Førsund, and Lychev (2012a), Roshdi, Van de Woestyne, Davtalab-Olyaie (2014) and Sueyoshi and Sekitani (2007a).

It is, however, worth noting that while the minimum face is a polytope in the VRS-based technology, it is an unbounded polyhedral cone in the CRS-based technology that is generated by the reference units in the GRS. Despite this structural difference between the two technologies[6], our results can still be successfully adapted for the case of CRS by removing the convexity constraint, i.e., $\sum_{j \in J} \lambda_j = 1$. This is because our approach is primarily based on finding a solution with the maximum number of positive components for a linear system of equations, and is independent of the existence of the convexity constraint, accordingly.

➢ *Dealing with negative input–output data*

Being independent of the data sets used, our proposed approach is free from the restricting assumption that the input–output data must be non-negative, which makes the identification of the GRS possible in the presence of negative data. From a practical point of view, this can be very beneficial since we deal with negative inputs and/or outputs in many empirical applications.

## 4. Illustration of the proposed approach

### 4.1. Numerical example

Let us consider a data set exhibited in Table 1 that consists of eight hypothetical DMUs with one input and one output. Based on these data, Fig. 1 depicts the frontier spanned in the two-dimensional input–output space.

**Table 1.** Input and output data for Example 4.1

|  | DMU$_1$ | DMU$_2$ | DMU$_3$ | DMU$_4$ | DMU$_5$ | DMU$_6$ | DMU$_7$ | DMU$_8$ |
|---|---|---|---|---|---|---|---|---|
| **Input** | 1 | 2 | 3 | 5 | 8 | 2 | 3 | 6 |
| **Output** | 2 | 5 | 6 | 8 | 8 | 1 | 3 | 4 |

To illustrate the application of our proposed approach, we first evaluate each DMU using the RAM model. Table 2 exhibits the efficiency score and the projection obtained for each DMU. The results reveal that DMU$_1$, DMU$_2$, DMU$_3$ and DMU$_4$ form the efficient frontier, and hence, are RAM-efficient. Amongst the inefficient DMUs (DMU$_5$, DMU$_6$, DMU$_7$ and DMU$_8$), DMU$_8$ has the minimum efficiency score of $\rho_8 = 0.643$.

---

[6] For more details about the facial structure of the CRS- and VRS-based technologies, see Davtalab-Olyaie et al. (2014), Davtalab-Olyaie, Roshdi, Partovi Nia, and Asgharian (2014) and Jahanshahloo, Roshdi, and Davtalab-Olyaie (2013).

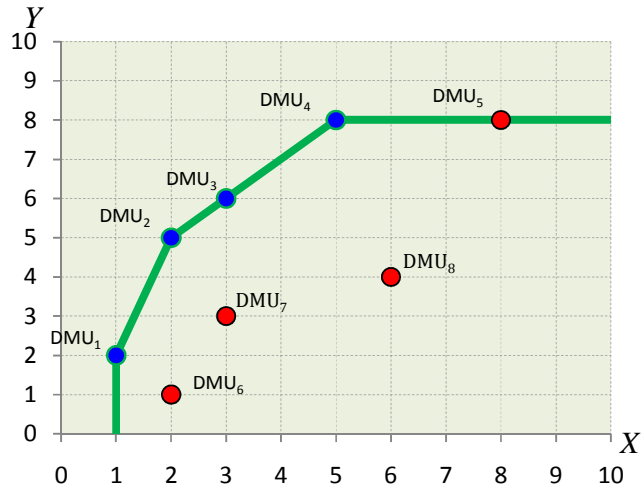

**Fig 1.** The production frontier

While DMU$_4$ and DMU$_2$ are found as the unique projections for DMU$_5$ and DMU$_6$, respectively; DMU$_7$ and DMU$_8$ do not have unique projections despite the fact that DMU$_3$ is determined as a projection for both. The projection sets of DMU$_7$ and DMU$_8$ are, respectively, the line segment connecting DMU$_2$ and DMU$_3$ and the line segment connecting DMU$_2$ and DMU$_4$ ( $R_1^- = 8$ and $R_1^+ = 8$ ). So, the minimum face associated with each of these units is the line segment connecting DMU$_2$ and DMU$_4$.

**Table 2.** The results for Example 4.1

| | | | DMU$_1$ | DMU$_2$ | DMU$_3$ | DMU$_4$ | DMU$_5$ | DMU$_6$ | DMU$_7$ | DMU$_8$ |
|---|---|---|---|---|---|---|---|---|---|---|
| **Model (3)** | | $\rho_o$ | 1 | 1 | 1 | 1 | 0.786 | 0.714 | 0.786 | 0.643 |
| | | $P$ | DMU$_1$ | DMU$_2$ | DMU$_3$ | DMU$_4$ | DMU$_4$ | DMU$_2$ | DMU$_3$ | DMU$_3$ |
| **Model (21)** | Ref. Weights | $\lambda_1'$ | 1 | | | | | | | |
| | | $\lambda_2'$ | | 1 | 0.500 | | | 1 | 0.500 | 0.333 |
| | | $\lambda_3'$ | | | 0.250 | | | | 0.250 | 0.333 |
| | | $\lambda_4'$ | | | 0.250 | 1 | 1 | | 0.250 | 0.333 |
| | | $P_o^*$ | DMU$_1$ | DMU$_2$ | DMU$_3$ | DMU$_4$ | DMU$_4$ | DMU$_2$ | DMU$_3$ | (3.333,6.333) |
| **RTS measurement** | | $u_o^*$ | -1 | 1.500 | 1 | 0.600 | 0.600 | 1.500 | 1 | 0.900 |
| | | $u_o^+ / u_o^-$ | -0.333 | 0 | 1 | 0.600 | 0.600 | 0 | 1 | 0.900 |
| | | RTS | IRS | CRS | DRS | DRS | DRS | CRS | DRS | DRS |

**CRS:** Constant Returns to Scale; **IRS:** Increasing Returns to Scale; **DRS:** Decreasing Returns to Scale.

Having obtained the efficiency score for each DMU, we use model (21) to identify its GRS and determine a relative interior point of its corresponding minimum face. The results are all presented in Table 2. Out of the four RAM-efficient DMUs, $DMU_1$, $DMU_2$ and $DMU_4$ are extreme-efficient and $DMU_3$ is non-extreme efficient (see Fig. 1).

Since $DMU_5$ and $DMU_6$ have unique projections (i.e., $DMU_4$ for $DMU_5$ and $DMU_2$ for $DMU_6$), the GRS for each of these units is exactly the same as its unique projection. Formally, $J_5^G = \{4\}$ and $J_6^G = \{2\}$.

Now, consider the case of $DMU_7$ suffering from the occurrence of problems Type I, II and III:

➤ Type I: The sets $\{DMU_3\}$, $\{DMU_2, DMU_4\}$ and $\{DMU_2, DMU_3, DMU_4\}$ are the three URSs for $DMU_7$ associated with the projection $DMU_3$. So, the MRS of $DMU_7$ associated with $DMU_3$ is $\{DMU_2, DMU_3, DMU_4\}$.

➤ Type II: $\Lambda_2$ is the line segment connecting $DMU_2$ and $DMU_3$.

➤ Type III resulted from Type II: The RTS of $DMU_7$ based upon the RTS of $DMU_2$ and $DMU_3$ (as its projections) is constant and decreasing, respectively.

As can be seen in Table 2, the GRS of $DMU_7$ consists of $DMU_2$, $DMU_3$ and $DMU_4$ with the respective weights of 0.5, 0.25 and 0.25, i.e., $J_7^G = \{2, 3, 4\}$. This finding confirms Corollary 3.1.1, i.e., the convex hull of $DMU_2$, $DMU_3$ and $DMU_4$ is $\Gamma_7^{\min}$. Moreover, $DMU_3$ is determined as a relative interior point of $\Gamma_7^{\min}$.

We now turn to measure the RTS of each inefficient DMU based upon the RTS of the relative interior point of its associated minimum face given in Table 2. To do so, we use the two–stage method of Banker et al. (2004) and examine the sign of the intercept of the supporting hyperplane(s) passing through the given relative interior point. In Stage 1, we solve the multiplier form of the BCC model and obtain the intercept; if the intercept is equal to zero, the RTS is constant. Otherwise, depending on the sign of the intercept, we solve an additional problem in Stage 2 to determine the RTS.

We apply the two–stage approach to each relative interior point and summarize the results in Table 2. As can be seen, out of the four inefficient units, the three units - $DMU_5$, $DMU_7$ and $DMU_8$ operate under decreasing RTS and one unit (i.e., $DMU_6$) under constant RTS.

### 4.2. Empirical application

To demonstrate the ready applicability of our proposed approach, we conduct an illustrative empirical analysis based on a real-life data set of 70 public schools in the United States, which was taken from Charnes et al. (1981). To carry out all the computations, we have developed a computer program using the GAMS optimization software.

The data consists of five inputs and three outputs. The inputs of schools are the education level of mother as measured in terms of percentage of high school graduates among female parents ($x_1$), the highest occupation of a family member according to a pre-arranged rating scale ($x_2$), the parental visit index representing the number of visits to the school site ($x_3$), the parent counseling index calculated from data on the time spent with child on school-related topics such as reading together, etc ($x_4$), and the number of teachers at a given site ($x_5$). The outputs are the Total Reading Score as measured by the Metropolitan Achievement Test ($y_1$), the Total Mathematics Score as measured by the Metropolitan Achievement Test ($y_2$), and the Coopersmith Self-Esteem Inventory, intended as a measure of self-esteem ($y_3$).

**Table 3.** Descriptive statistics of efficiency scores and projections obtained by the RAM model for inefficient schools

| DMU | $\rho_o$ | $\hat{x}_1$ | $\hat{x}_2$ | $\hat{x}_3$ | $\hat{x}_4$ | $\hat{x}_5$ | $\hat{y}_1$ | $\hat{y}_2$ | $\hat{y}_3$ |
|------|------|--------|--------|--------|--------|--------|--------|--------|--------|
| Min  | 0.844 | 4.379  | 2.233  | 7.958  | 8.170  | 2.274  | 8.817  | 9.369  | 6.350  |
| Max  | 0.985 | 39.969 | 14.650 | 51.902 | 52.132 | 8.314  | 54.530 | 63.557 | 39.100 |
| Mean | 0.938 | 15.634 | 8.396  | 28.938 | 29.847 | 4.947  | 30.319 | 36.529 | 22.367 |

In order to apply model (21), we have to obtain the efficiency score of $\rho_o$ beforehand. Thus, we first make the efficiency assessment of the schools using the RAM model. The results show that 26 (37%) schools are RAM-efficient. For the remaining 44 (63%) inefficient schools, Table 3 provides summary statistics of the efficiency scores and the inputs and outputs of the projection points. It can be observed that the values of $\rho_o$ range from 0.844 to 0.965 with the mean efficiency score of 0.938.

**Table 4.** The GRSs and RTSs of the inefficient schools

| DMU | Global Reference Set | | | | | Returns to Scale | | |
|---|---|---|---|---|---|---|---|---|
| S1 | S52( 0.652) | S58( 0.171) | S59( 0.177) | | | 0.075 | 0.003 | DRS |
| S2 | S44( 0.250) | S58( 0.750) | | | | 0.004 | 0 | CRS |
| S3 | S44( 0.201) | S52( 0.275) | S58( 0.525) | | | 0.004 | 0 | CRS |
| S4 | S44( 0.036) | S58( 0.954) | S59( 0.010) | | | 0.006 | 0 | CRS |
| S6 | S58( 0.862) | S62( 0.138) | | | | -0.006 | 0 | CRS |
| S7 | S58( 0.944) | S59( 0.056) | | | | 0.044 | 0 | CRS |
| S8 | S44( 0.003) | S58( 0.739) | S59( 0.258) | | | 0.024 | 0 | CRS |
| S9 | S58( 0.831) | S59( 0.169) | | | | -0.003 | 0 | CRS |
| S10 | S44( 0.159) | S52( 0.555) | S58( 0.236) | S59( 0.049) | | 0.023 | 0.003 | DRS |
| S13 | S58( 0.828) | S59( 0.172) | | | | 0.121 | 0 | CRS |
| S14 | S58( 0.104) | S62( 0.624) | S69( 0.272) | | | -0.014 | 0 | CRS |
| S16 | S44( 0.337) | S58( 0.558) | S59( 0.105) | | | 0.003 | 0 | CRS |
| S19 | S44( 0.397) | S58( 0.560) | S59( 0.043) | | | 0.003 | 0 | CRS |
| S23 | S44( 0.478) | S58( 0.500) | S59( 0.022) | | | 0.021 | 0 | CRS |
| S25 | S44( 0.234) | S58( 0.299) | S62( 0.467) | | | -0.047 | 0 | CRS |
| S26 | S44( 0.201) | S58( 0.698) | S59( 0.101) | | | 0.05 | 0 | CRS |
| S28 | S20( 0.047) | S20( 0.047) | S27( 0.369) | S47( 0.084) | S62( 0.104) | -0.096 | -0.061 | IRS |
| S29 | S62( 0.627) | S69( 0.373) | | | | -0.039 | 0 | CRS |
| S30 | S58( 0.973) | S59( 0.027) | | | | 0.05 | 0 | CRS |
| S31 | S58( 0.990) | S59( 0.010) | | | | 0.004 | 0 | CRS |
| S33 | S44( 0.288) | S58( 0.563) | S59( 0.150) | | | 0.021 | 0 | CRS |
| S34 | S58( 0.768) | S59( 0.232) | | | | -0.001 | 0 | CRS |
| S36 | S58( 0.997) | S59( 0.003) | | | | 0.109 | 0 | CRS |
| S37 | S58( 0.978) | S59( 0.022) | | | | 0.001 | 0 | CRS |
| S39 | S20( 0.396) | S27( 0.461) | S44( 0.143) | | | 0.058 | 0 | CRS |
| S40 | S52( 0.239) | S58( 0.666) | S69( 0.095) | | | 0.012 | 0 | CRS |
| S41 | S44( 0.229) | S58( 0.767) | S59( 0.004) | | | 0.034 | 0 | CRS |
| S42 | S17( 0.482) | S20( 0.482) | S58( 0.036) | | | 0.01 | 0 | CRS |
| S43 | S44( 0.114) | S58( 0.776) | S59( 0.110) | | | 0.029 | 0 | CRS |
| S46 | S44( 0.032) | S58( 0.710) | S59( 0.259) | | | 0.023 | 0 | CRS |
| S50 | S52( 0.545) | S58( 0.399) | S59( 0.056) | | | 0.037 | 0.004 | DRS |
| S51 | S58( 0.274) | S62( 0.726) | | | | 0.017 | 0 | CRS |
| S53 | S58( 0.948) | S59( 0.052) | | | | 0.045 | 0 | CRS |
| S55 | S15( 0.010) | S44( 0.117) | S52( 0.076) | S58( 0.554) | S69( 0.243) | -0.001 | 0 | CRS |
| S57 | S52( 0.431) | S58( 0.499) | S69( 0.070) | | | 0.017 | 0 | CRS |
| S60 | S44( 0.091) | S58( 0.143) | S62( 0.766) | | | -0.011 | 0 | CRS |
| S61 | S58( 0.009) | S62( 0.447) | S69( 0.544) | | | 0.025 | 0 | CRS |
| S63 | S44( 0.017) | S52( 0.012) | S58( 0.437) | S62( 0.058) | S69( 0.476) | 0.027 | 0 | CRS |
| S64 | S44( 0.032) | S58( 0.863) | S59( 0.105) | | | 0.012 | 0 | CRS |
| S65 | S44( 0.006) | S52( 0.075) | S58( 0.219) | S62( 0.224) | S69( 0.476) | 0.016 | 0 | CRS |
| S66 | S44( 0.052) | S58( 0.915) | S59( 0.033) | | | 0.042 | 0 | CRS |
| S67 | S44( 0.063) | S58( 0.834) | S59( 0.103) | | | -0.006 | 0 | CRS |
| S68 | S44( 0.147) | S58( 0.249) | S59( 0.134) | S62( 0.471) | | -0.004 | 0 | CRS |

| S70 | S44( 0.027) | S58( 0.640) | S62( 0.333) | | 0.053 | 0 | CRS |

As is well known, the RAM model is not by itself able to find out all the possible reference units for inefficient DMUs. Therefore, we apply model (21) for each inefficient school to identify the reference schools in its GRS as appropriate benchmarks. For the 44 inefficient schools, Table 4 shows all the reference schools together with their corresponding weights. For example, the efficient schools S44, S58 and S59 with the respective weights of 0.052, 0.915 and 0.033 appear in the GRS of the most inefficient school S66. This means that S66's target inputs and outputs are a linear combination of S44, S58 and S59's inputs and outputs. Thus, in order for S66 to become efficient, it must adjust its inputs and outputs so that it produces $0.052 \times \mathbf{y}_{S44} + 0.915 \times \mathbf{y}_{S58} + 0.033 \times \mathbf{y}_{S59}$ amount of output by consuming $0.052 \times \mathbf{x}_{S44} + 0.915 \times \mathbf{x}_{S58} + 0.033 \times \mathbf{x}_{S59}$ amount of input.

Now, we proceed to measure the RTS of the inefficient schools. For each inefficient DMU, we first use (24) to obtain a relative interior point of its corresponding minimum face. The statistics of the results are given in Table 4. Applying the two–stage method of Banker et al. (2004) to the obtained interior points of Table 5 yields the results that are reported in the last three columns of Table 4. As the results show, among 44 inefficient branches, 3 (7%) schools have decreasing RTS status, 1 (2%) has increasing RTS status, and the remaining 40 (91%) schools have constant RTS status.

**Table 5.** Descriptive statistics of relative interior points of the minimum faces

| DMU | $\hat{x}_1$ | $\hat{x}_2$ | $\hat{x}_3$ | $\hat{x}_4$ | $\hat{x}_5$ | $\hat{y}_1$ | $\hat{y}_2$ | $\hat{y}_3$ |
|------|--------|--------|--------|--------|-------|--------|--------|--------|
| Min | 4.379 | 2.233 | 7.958 | 8.170 | 2.274 | 8.817 | 9.369 | 6.350 |
| Max | 39.969 | 14.650 | 51.902 | 52.132 | 8.314 | 54.530 | 63.557 | 39.100 |
| Mean | 15.634 | 8.396 | 28.938 | 29.847 | 4.947 | 30.319 | 36.529 | 22.367 |

## 5. Summary and concluding remarks

The current study is mainly concerned with the identification of all the possible reference units of an evaluated inefficient DMU. It is also interested in its application for the measurement of RTS in non-radial DEA models. Corresponding to a given projection of the DMU under evaluation, first, two basic notions were introduced: i) URS: the set of efficient DMUs that are active in a specific convex combination generating this projection, and ii) MRS: the union of all the URSs associated with this projection. Then, the notion of GRS was defined as the union of the MRSs associated with all projections of the evaluated DMU. With the help of the introduced notions, it was demonstrated that the convex hull of the GRS is equal to the minimum face, from which it was immediately concluded that the minimum face is a polytope.

Three types of multipleness may occur in any non-radial DEA model: multiple URSs (Type I), multiple projections (Type II), and multiple supporting hyperplanes (Type III). The occurrences of problems Type I and II cause difficulties in the identification of all the possible reference DMUs. The difficulty in the measurement of RTS arises mainly from problem Type III, which itself originates from two sources: problem Type II and the non-full dimensionality of the minimum face.

To deal effectively with the simultaneous occurrence of problems Type I and II, an LP-based approach was proposed to identify the GRS. Our proposed approach has several advantages over the existing ones. First, since it requires the execution of a single LP problem, it is computationally more efficient than the existing ones for easy implementation in practical applications. Second, using the simplex algorithm adopted for solving the LP problems dealing with upper-bounded variables, the computational efficiency of our approach can be substantially improved. Third, as our proposed approach is primal-based, its computational efficiency is higher than that of the previous primal–dual methods. Fourth, our proposed approach is more general in the sense that it can easily be applied to both radial and non-radial DEA models (e.g., the BCC and the additive models).

To estimate the RTS, a method in the non-radial DEA framework is also proposed to deal effectively with problem Type III, which arises either from problem Type II or from the non-full dimensionality of the minimum face. A key outcome of the LP problem proposed to identify the GRS is that it generates a projection in the relative interior point of the minimum face. Using this projection for determining the RTS of the evaluated inefficient DMU, our proposed method overcomes the difficulty arising from problem Type II. This is because each supporting hyperplane binding at the used projection is determined through the GRS, but not through a specific MRS. To cope with the difficulty arising from the non-full dimensionality of the minimum face, our proposed method employs the indirect method of Banker et al. (2004) or the direct method of Førsund et al. (2007).

## Acknowledgments

We wish to thank Robert G. Dyson (Editor), three anonymous referees of the Journal and Abbas Valadkhani, whose invaluable inputs and comments considerably improved an earlier version of this article. We are also thankful to Arash Moradi for his careful proofread of the article. The usual caveat applies.

## Appendix A

**Proof of Lemma 3.1.1** Let $(\bar{\mathbf{x}}, \bar{\mathbf{y}})$ be a relative interior point of $\Lambda_o$. Since $(\bar{\mathbf{x}}, \bar{\mathbf{y}})$ is RAM-efficient, the SCSCs of linear programming implies the existence of a strong supporting hyperplane of $T_{VRS}^{DEA}$ that passes through this point. Without loss of generality, let $H^S$ be such a hyperplane whose associated strong face $F^S := H^S \cap T_{VRS}^{DEA}$ is of minimum dimension. By Theorem 6.4 in Rockafellar (1970), the convexity of $\Lambda_o$ infers that $H^S$ is binding at all the DMUs in $R_o^G$. Therefore, $conv\left(R_o^G\right) \subseteq H^S$, indicating that $conv\left(R_o^G\right) \subseteq F^S$.

According to Theorem 2 in Davtalab-Olyaie et al. (2014), $F^S$ is a polytope (bounded polyhedral set). Hence, the equality holds if $\Lambda_o \cap ri\left(F^S\right) \neq \varnothing$, since this relation implies that all the observed DMUs on $F^S$ belong to $R_o^G$. Assume, by the way of contradiction, that $\Lambda_o \cap ri\left(F^S\right) = \varnothing$ or, equivalently, $\Lambda_o \subseteq \partial\left(F^S\right)$. Then, there exists a unique face of $F^S$ of minimum dimension, namely $K^S$, for which $\Lambda_o \subseteq K^S \subseteq \partial\left(F^S\right) \subsetneqq F^S$. The face $K^S$ is a strong face of $T_{VRS}^{DEA}$ (Rockafellar, 1970), whose dimension is less than that of $F^S$. This is a contradiction and, thus, the proof is complete by Definition 2.1.1. $\qquad\qquad\square$

**Proof of Theorem 3.1.1** As proved in Lemma 3.1.1, $conv\left(R_o^G\right)$ is a strong face of $T_{VRS}^{DEA}$ containing $\Lambda_o$. Thus, according to the definition of $\Gamma_o^{\min}$, it will suffice to show that $conv\left(R_o^G\right) \subseteq \Gamma_o^{\min}$. By the definition of a face, there exists a supporting hyperplane, namely $H^{\min}$, such that $\Gamma_o^{\min} = H^{\min} \cap T_{VRS}^{DEA}$. Since $H^{\min}$ is binding at each projection, it passes through each DMU in $R_o^G$. Then, the convexity of $H^{\min}$ follows that $conv\left(R_o^G\right) \subseteq H^{\min}$, which completes the proof. $\qquad\qquad\square$

**Proof of Theorem 3.1.2** By the assumption, we have $\left(\boldsymbol{\lambda}',\mathbf{s}^{-\prime},\mathbf{s}^{+\prime}\right) \in \Omega_o$. Hence, from (13), we need only to prove that $R_o^G \subseteq \left\{j \big| \lambda_j' > 0\right\}$, which is equivalent to demonstrating that $\boldsymbol{\lambda}'$ takes positive values in any positive component of $\boldsymbol{\lambda}$ in any $\left(\boldsymbol{\lambda},\mathbf{s}^{-},\mathbf{s}^{+}\right) \in \Omega_o$.

By the way of contradiction, assume that there exists an element $\left(\boldsymbol{\lambda}'',\mathbf{s}^{-\prime\prime},\mathbf{s}^{+\prime\prime}\right) \in \Omega_o$ and an index $j_h \in \left\{j \big| \lambda_j' > 0\right\}$ for which $\lambda_{j_h}' = 0$. Then, let $\left(\hat{\boldsymbol{\lambda}},\hat{\mathbf{s}}^{-},\hat{\mathbf{s}}^{+}\right)$ be a strict convex combination of the elements $\left(\boldsymbol{\lambda}',\mathbf{s}^{-\prime},\mathbf{s}^{+\prime}\right)$ and $\left(\boldsymbol{\lambda}'',\mathbf{s}^{-\prime\prime},\mathbf{s}^{+\prime\prime}\right)$. Since $\Omega_o$ is convex, $\left(\hat{\boldsymbol{\lambda}},\hat{\mathbf{s}}^{-},\hat{\mathbf{s}}^{+}\right) \in \Omega_o$ and $\left\{j \big| \hat{\lambda}_j > 0\right\} = \left\{j \big| \lambda_j' > 0\right\} \cup \left\{j \big| \lambda_j'' > 0\right\}$. Consequently, $\left\{j \big| \lambda_j' > 0\right\} \subsetneqq \left\{j \big| \hat{\lambda}_j > 0\right\}$, which contradicts the assumption that $\boldsymbol{\lambda}'$ has the maximum number of positive components. Thus, the proof is complete. $\quad\square$

**Proof of Lemma 3.2.1** Let $\left(\mathbf{u}^{*}+\mathbf{w}^{*},\mathbf{v}^{*}\right)$ be an optimal solution to (18). Then $\left(\mathbf{u}^{*}+\mathbf{w}^{*},\mathbf{v}^{*}\right) \in X$. Without loss of generality, assume that $\mathbf{u}^{*}+\mathbf{w}^{*} = \left(u_1^{*}+w_1^{*},...,u_k^{*}+w_k^{*},0,...,0\right)$. We claim that the number of positive components of $\mathbf{u}^{*}+\mathbf{w}^{*}$ is maximum or, equivalently, that $\mathbf{u}^{*}+\mathbf{w}^{*}$ takes positive values in any positive component of $\mathbf{u}$ in any feasible solution $\left(\mathbf{u},\mathbf{v}\right) \in X$ because of the convexity of $X$. To prove our claim, assume on the contrary that there exists an element of $X$, namely $\left(\hat{\mathbf{u}},\hat{\mathbf{v}}\right)$, and a set of indices, namely $\left\{k+1,...,l\right\}$, such that $\hat{u}_{k+1} > 0,...,\hat{u}_l > 0$. Since the system of equations in (18) is homogeneous, without loss of generality, it can be assumed that $\hat{u}_j \le 1$, $j = 1,...,q_1$.

Since $\left(\mathbf{u}^{*}+\mathbf{w}^{*},\mathbf{v}^{*}\right) \in X$ and $\left(\hat{\mathbf{u}},\hat{\mathbf{v}}\right) \in X$, we have $\left(\mathbf{u}^{*}+\mathbf{w}^{*}+\hat{\mathbf{u}},\mathbf{v}^{*}+\hat{\mathbf{v}}\right) \in X$ thereby

$$\left[\sum_{j=1}^{k}\mathbf{a}_j\left(\left(u_j^{*}+\hat{u}_j\right)+w_j^{*}\right)+\sum_{j=k+1}^{q_1}\mathbf{a}_j\hat{u}_j\right]+\sum_{j=1}^{q_2}\mathbf{b}_j\left(v_j^{*}+\hat{v}_j\right)=\mathbf{0}. \tag{A.1}$$

Based on (A.1), we define

$$u_j' = \begin{cases} u_j^{*}+\hat{u}_j, & j=1,...,k, \\ 0, & j=k+1,..,q_1; \end{cases} \quad w_j' = \begin{cases} w_j^{*}, & j=1,...,k, \\ \hat{u}_j, & j=k+1,..,q_1; \end{cases} \quad v_j' = v_j^{*}+\hat{v}_j, \ j=1,...,q_2. \tag{A.2}$$

Then, $\left(\mathbf{u}', \mathbf{w}', \mathbf{v}'\right)$ is a feasible solution of (18) with the objective value of $\sum_{j=1}^{q_1} w_j{}'$, which is strictly greater than $\sum_{j=1}^{q_1} w_j^*$. This contradicts the optimality of $\left(\mathbf{u}^* + \mathbf{w}^*, \mathbf{v}^*\right)$, and hence proves our claim. $\qquad\square$

**Proof of Theorem 3.2.1** Let $\left(u_j^* + w_j^*, j = 1, ..., q_1 + 1, v_j^*, j = 1, ..., q_2\right)$ be an optimal solution to (20). By Lemma 3.2.1, this solution has the maximum number of positive components. Therefore, $u_{q_1+1}^* + w_{q_1+1}^*$ is positive as a result of the assumption that $X \neq \varnothing$. Consequently, the first constraint of (20) at optimality can be equivalently rewritten as

$$\sum_{j=1}^{q_1} \mathbf{a}_j \left( \frac{u_j^* + w_j^*}{u_{q_1+1}^* + w_{q_1+1}^*} \right) + \sum_{j=1}^{q_2} \mathbf{b}_j \left( \frac{v_j^*}{u_{q_1+1}^* + w_{q_1+1}^*} \right) = \mathbf{d} \,, \qquad (A.3)$$

which completes the proof. $\qquad\square$